\documentclass[a4paper]{amsart}
\usepackage{amssymb,amscd}
\usepackage{fullpage}

\newtheorem{thm}{Theorem}[section]
\newtheorem{prop}[thm]{Proposition}
\newtheorem{lem}[thm]{{Lemma}}
\newtheorem{defn}[thm]{Definition}

\numberwithin{equation}{section}

\def\.{\cdot}
\def\<{\left\langle}
\def\>{\right\rangle}
\def\({\left(}
\def\){\right)}

\renewcommand{\leq}{\leqslant}
\renewcommand{\geq}{\geqslant}
\renewcommand{\phi}{\varphi}

\def\bar#1{\overline{#1}}

\def\H{\mathbb H}
\def\L{\mathbb\L}

\def\subset{\subseteq}

\def\epsilon{\varepsilon}

\def\index{\operatorname{index}}

\def\C{\mathbb C}

\def\Z{\mathbb Z}

\def\H{\mathbb H}

\def\S{\mathbb S}
\def\ds{\displaystyle}

\def\virt-dim{\operatorname{virt-dim}}

\def\rs{\operatorname{res}}
\def\Id{\operatorname{Id}}

\def\e{\operatorname{e}}

\def\exp{\operatorname{exp}}

\def\dim{\operatorname{dim}}

\title[Complex Cobordism of Hilbert Manifolds]
{On smooth Chas-Sullivan loop product in Quillen's Geometric
Complex Cobordism of Hilbert Manifolds}
\author{Cenap \"{O}zel}
\address{AIBU Golkoy Kampusu, Bolu 14280, Turkey.}
\email{cenap@ibu.edu.tr}
\subjclass{Algebraic topology,Global Analysis}
\date{September 02,2003}

\begin{document}
\begin{abstract}
In \cite{baker-ozel}, by using Fredholm index we developed a
version of Quillen's geometric cobordism theory for infinite
dimensional Hilbert manifolds. This cobordism theory has a graded
group structure under topological union operation and has
push-forward maps for complex orientable Fredholm maps. In
\cite{cenap-isr}, by using Quinn's Transversality Theorem
\cite{Quinn}, it has been  shown that this cobordism theory has a
graded ring structure under transversal intersection operation and
has pull-back maps for smooth maps. It has been shown that the
Thom isomorphism in this theory was  satisfied for finite
dimensional vector bundles over separable Hilbert manifolds and
the projection formula for Gysin maps has been proved. In
\cite{chas}, Chas and Sullivan described an intersection product
on the homology of loop space $LM$. In \cite{cohen}, R. Cohen and
J. Jones described a realization of the Chas-Sullivan loop product
in terms of a ring spectrum structure on the Thom spectrum of a
certain virtual bundle over the loop space. In this paper, we will
extend this product on cobordism and bordism theories.
\end{abstract}
\maketitle

\section{The Fredholm Index and Complex Cobordism of Hilbert Manifolds.}
In \cite{Quillen}, Quillen gave a geometric interpretation of
cobordism groups which suggests a way of defining the cobordism of
separable Hilbert manifolds equipped with suitable structure. In
order that such a definition be sensible, it ought to reduce to
his for finite dimensional manifolds and smooth maps of manifolds
and be capable of supporting reasonable calculations for important
types of infinite dimensional manifolds such as homogeneous spaces
of free loop groups of finite dimensional Lie groups.

\subsection{Cobordism of separable Hilbert  manifolds.} By a
manifold, we mean a smooth manifold modelled on a separable
Hilbert space; see Lang \cite{lang} for details on infinite
dimensional manifolds. The facts about Fredholm map can be found
in \cite{conway}.
\begin{defn}\label{zirto}
Suppose that $f: X\rightarrow Y$ is a proper Fredholm map with
even index at each point. Then $f$ is an \emph{admissible complex
orientable map} if there is a smooth factorization
$$f: X\stackrel{\tilde{f}}\rightarrow \xi \stackrel{q}\rightarrow Y,$$
where $q:\xi \rightarrow Y$ is a finite dimensional smooth complex
vector bundle and $\tilde{f}$ is a smooth embedding endowed with a
complex structure on its normal bundle $\nu (\tilde{f})$.
\par
A complex orientation for a Fredholm  map $f$ of odd index is
defined to be one for the map $ (f,\epsilon): X\rightarrow Y\times
\mathbb R $ given by  $(f,\epsilon) (x) = (f(x),0)$ for every $x
\in X$. At $x\in  X$, $\index (f, \epsilon)_x = (\index f_x) - 1$.
Also the finite dimensional  complex vector bundle $\xi$ in the
smooth factorization will be replaced by $\xi \times \mathbb R$.
\end{defn}
Suppose that $f$ is an admissible complex orientable map. Then
since the map $f$ is the Fredholm and $\xi$ is a finite
dimensional vector bundle, we see $\tilde{f}$ is also a Fredholm
map. By the surjectivity of $q$,
$$\index \tilde{f} = \index f - \dim \xi.$$
Before we give a notion of equivalence of such factorizations
$\tilde{f}$ of $f$, we want to give some definitions.
\begin{defn}
Let $X$, $Y$ be the smooth separable Hilbert manifolds and $F: X
\times \mathbb R \rightarrow Y$ a smooth map. Then we will say
that $F$ is an \emph{isotopy} if it satisfies the following
conditions.
\begin{enumerate}
\item For every $t \in \mathbb R$, the map $F_t$ given by
$F_t (x) = F(x, t)$ is an embedding.
\item There exist numbers
$t_0 < t_1$ such that $F_t = F_{t_0}$ for all $t \leq t_0$ and
$F_t = F_{t_1}$ for all $t \geq t_1$.
\end{enumerate}
The closed interval $[t_0, t_1]$ is called a \emph{proper domain}
for the isotopy. We say that two embeddings $f: X \rightarrow Y$
and $g: X \rightarrow Y$ are \emph{isotopic} if there exists an
isotopy $F_t: X\times \mathbb R \rightarrow Y$ with proper domain
$[t_0, t_1]$ such that $f = F_{t_0}$ and $g = F_{t_1}$.
\end{defn}
\begin{prop}(see \cite{lang})\label{lang}
The relation of isotopy between smooth embeddings is an
equivalence relation.
\end{prop}
\begin{defn}
Two factorizations $f: X\stackrel{\tilde{f}}\rightarrow \xi
\stackrel{q}\rightarrow Y$ and $f:
X\stackrel{\tilde{f'}}\rightarrow \xi' \stackrel{q'}\rightarrow Y$
are \emph{equivalent} if $\xi$ and $\xi'$ can be embedded as
subvector bundles of a vector bundle $\xi''\rightarrow Y$ such
that $\tilde{f}$ and $\tilde{f'}$ are isotopic in $\xi''$ and this
isotopy is compatible with the complex structure on the normal
bundle. That is, there is an isotopy $F$ such that for all $t \in
[t_0, t_1]$, $F_t:  X \rightarrow \xi''$ is endowed with a complex
structure on its normal bundle which matches that of $\tilde{f}$
and $\tilde{f'}$ in $\xi''$ at $t_0$ and $t_1$ respectively.
\end{defn}
By Proposition \ref{lang}, we have
\begin{prop}
The relation of equivalence of admissible complex orientability of
proper Fredholm maps between separable Hilbert manifolds is an
equivalence relation.
\end{prop}
This generalizes  Quillen's notion of complex orientability for
maps of finite dimensional manifolds. We can also define a notion
of cobordism of admissible complex orientable maps between
separable Hilbert manifolds. First we recall some ideas on the
transversality.

\begin{defn}Let $f_1: M_1 \rightarrow N, f_2: M_2 \rightarrow N$
be smooth maps between Hilbert manifolds. Then $f_1$ and $f_2$ are
\emph{transverse} at $y \in N$ if
\[
df_1 (T_{x_1} M_1) + df_2 (T_{x_2}M_2) = T_y N
\]
whenever $f_1 (x_1) = f_2(x_2) = y$. The maps $f_1$ and $f_2$ are
said to be \emph{transverse} if they are transverse at every point
of $N$.
\end{defn}
\begin{lem}
Smooth maps $f_i: M_i \rightarrow N (i = 1,2)$ are transverse if
and only if $f_1\times f_2 : M_1 \times M_2 \rightarrow N\times N$
is transverse to the diagonal map $\Delta :  N\rightarrow N\times
N$.
\end{lem}
\begin{defn}
Let $f_1: M_1 \rightarrow N, f_2: M_2 \rightarrow N$ be transverse
smooth maps between smooth Hilbert manifolds. The
\emph{topological pullback}
\[
M_1 \prod_N M_2 = \{ (x_1 , x_2) \in M_1 \times M_2 : f_1(x_1) =
f_2(x_2)\}
\]
is a submanifold of $M_1 \times M_2$ and the diagram
\[
\begin{CD}
M_1 \prod_N M_2  @>{f_2}^*(f_1)>> M_2  \\
@VV{f_1}^*(f_2)V         @VVf_2V \\
M_1 @>f_1>> N
\end{CD}
\]
is commutative, where the map ${f_i}^* (f_j)$ is pull-back of
$f_j$ by $f_i$.
\end{defn}

\begin{defn}
Let $f_i: X_i \rightarrow  Y (i= 0,1)$ be admissible complex
oriented maps. Then $f_0$ is \emph{cobordant} to $f_1$ if there is
an admissible complex orientable map $h: W \rightarrow Y\times
\mathbb R$ such that the maps $\epsilon_i: Y\rightarrow Y \times
\mathbb R \,\text{given by} \,\epsilon_i (y) = (y,i)$ for $i= 0, 1
$, are transverse to $h$ and the pull-back map ${\epsilon_i}^* h$
is equivalent to $f_i$. The cobordism class of $f: X\rightarrow Y$
will be denoted by $[X, f]$.
\end{defn}
\begin{prop}\label{grem}
If $f: X\rightarrow Y$ is an admissible complex orientable map and
$g: Z\rightarrow Y$ a smooth map transverse to $f$, then the
pull-back map
\[
{g^*}(f):  Z\prod_Y X\rightarrow Z
\]
is an admissible complex orientable map with finite dimensional
pull-back vector bundle
\[
g^*(\xi) = Z\prod_Y \xi = \{(z,v)\in Z\times \xi : g(z) = q(v)\}
\]
in the factorization of $g^* (f)$, where $q: \xi\rightarrow Y$ is
the finite-dimensional complex vector bundle in the factorization
of $f$ as in Definition \ref{zirto}.
\end{prop}
The next result was proved in \cite{ozel} by essentially the same
argument as in the finite dimensional situation using the Implicit
Function Theorem \cite{lang}.
\begin{thm}
Cobordism is an equivalence relation.
\end{thm}

\begin{defn}
For a separable Hilbert manifold $Y$, $\mathcal U^d (Y)$ is the
set of cobordism classes of  the admissible complex orientable
proper Fredholm maps of index $-d$.
\end{defn}
My next result is the following.

\begin{thm}\label{iain}
If $f: X\rightarrow Y$ is an admissible complex orientable
Fredholm map of index $d_1$ and $g: Y\rightarrow Z$ is an
admissible complex orientable Fredholm map of index $d_2$, then
$g\circ f: X\rightarrow Z$ is an admissible complex orientable
Fredholm map with index $d_1 + d_2$.
\end{thm}

Let $g:Y \rightarrow Z$ be an admissible complex orientable
Fredholm map of index $r$. By Theorem \ref{iain}, we have
\emph{push-forward, or Gysin map}
\[
g_{*}: \mathcal U^d (Y) \rightarrow \mathcal U^{d + r} (Z)
\]
given by $g_{*} ([X, f]) = ([X, g \circ f])$.

We show in \cite{ozel} that it is well-defined. If $g':
Y\rightarrow Z$ is a second map cobordant to $g$ then ${g'}_* =
g_*$; in particular, if $g$ and $g'$ are homotopic through proper
Fredholm maps they induce the same Gysin maps. Clearly, we have
$(h\circ g)_{*} = h_{*}g_{*}$ for admissible complex orientable
Fredholm maps $h,g$ and $\Id_{*} = \Id$.

The graded cobordism set $\mathcal U^* (Y)$ of the separable
Hilbert manifold $Y$ has a group structure given as follows. Let
$[X_1 , f_1]$ and $[X_2 , f_2]$ be cobordism classes. Then $[X_1 ,
f_1] + [X_2 , f_2]$ is the class of the map $f_1 \sqcup f_2: X_1
\sqcup X_2 \rightarrow Y$, where $X_1 \sqcup X_2$ is the
topological sum (disjoint union) of $X_1$ and $X_2$.  We show in
\cite{ozel} that this sum is well-defined. As usual, the class of
the empty set $\emptyset$ is the zero  element of  the cobordism
set and the negative of $[X,f]$ is itself with the opposite
orientation on the normal bundle of the embedding $\tilde{f}$.
Then we have
\begin{thm}
The graded cobordism set $\mathcal U^* (Y)$ of the admissible
complex orientable maps of $Y$ is a graded abelian group.
\end{thm}
Now we define relative cobordism .
\begin{defn}
If $A$ is a finite dimensional submanifold of $Y$, the relative
cobordism set $\mathcal U^* (Y,A)$ is the set of the admissible
complex orientable maps of $Y$ whose images lie in $Y-A$.
\end{defn}
More generally,
\begin{thm}\label{relative}
Let $A$ be a finite dimensional submanifold of $Y$. Then the
relative cobordism set $\mathcal U^* (Y,A)$ is a graded abelian
group and there is a homomorphism $\kappa^* :\mathcal U^*
(Y,A)\rightarrow \mathcal U^* (Y)$ by
$\kappa^*[M\stackrel{h}\rightarrow Y] =[M\stackrel{h}\rightarrow
Y]$ with $h(M)\subset Y-A$.
\end{thm}

If our cobordism  functor $\mathcal U^*( \,)$ of admissible
complex orientable Fredholm maps is restricted to finite
dimensional Hilbert manifolds, it agrees Quillen's complex
cobordism functor $MU^* (\,)$.
\begin{thm}
For finite dimensional separable Hilbert manifolds $A\subset Y$,
there is a natural isomorphism
$$
\mathcal U^*(Y,A)\cong MU^*(Y,A).
$$
\end{thm}
\subsection{Transversal approximations, contravariance and cup products.} We
would like to define a product structure on the graded cobordism
group $\mathcal U^* (Y)$. Given  cobordism classes $[X_1 , f_1]
\in \mathcal U^{d_1}(Y_1)$ and $[X_2 , f_2]\in \mathcal
U^{d_2}(Y_2)$, their external product is
\[
[X_1 , f_1]\times [X_2 , f_2 ] = [X_1 \times X_2 ,f_1 \times f_2]
\in \mathcal U^{d_1 + d_2} (Y_1\times Y_2).
\]
Although there is the external product in the category of
cobordism of separable Hilbert manifolds, we can not necessarily
define an internal product on $U^{*}(Y)$ unless $Y$ is a finite
dimensional manifold. However,if admissible complex orientable
Fredholm map $f_1 \times f_2:X_1 \times X_2\rightarrow Y\times Y$ is transverse to
the diagonal imbedding $\Delta:Y\rightarrow Y\times Y$, then we do have an
internal (cup) product
\[
[X_1 , f_1] \cup [X_2 , f_2] = \Delta^* [X_1\times X_2, f_1\times
f_2].
\]
If $Y$ is finite dimensional, then by Haefliger and Thom's Transversality
Theorem in \cite{thom}, every complex orientable map to $Y$ has a
transverse approximation, hence the cup product $\cup$ induces a
graded ring structure on $\mathcal U^* (Y)$. The unit element $1$
is represented by the identity map $Y\rightarrow Y$ with index
$0$. However F. Quinn \cite{Quinn} proved the generalization of
Thom's Transversality Theorem for separable Hilbert manifolds
using smooth transversal approximations of Sard functions in fine
topology.

\par
By Quinn's Transversality Theorem,
a smooth map (even continuous map)
$g: Z\rightarrow Y$ can be deformed to a smooth map $g':
Z\rightarrow Y$ by a small correction until it is transverse to an
admissible complex orientable map $f: X \rightarrow Y$. It is
obvious that they are homotopic each other. By definition of
Cobordism and Proposition \ref{grem}, the cobordism functor is
contravariant for any smooth map between separable Hilbert
manifolds.
\begin{thm}\label{son}
Let $f: X\rightarrow Y$ be an admissible complex oriented map and
let $g: Z\rightarrow Y$ be a smooth (may be continuous) map. Then
the cobordism class of the pull-back $Z \prod_Y X \rightarrow Z$
depends only on the cobordism class of $f$, hence there is a map
$g^*: \mathcal U^d (Y)\rightarrow \mathcal U^d (Z)$ given by
$$g^* [X, f] =g'^* [X, f] = [Z\prod_Y X, {g'}^*(f)],$$ where $g'$ is
a smooth $\varepsilon$-approximation of g which is transverse to
$f$. Moreover, $g^*$ depends only on the homotopy class of $g$.
\end{thm}

\par

Let turn back the interior(cup) products in $\mathcal U^*$. Given
cobordism classes $[X_1,f_1] \in \mathcal U^{d_1}(Y_1)$ and $[X_2
, f_2]\in \mathcal U^{d_2}(Y_2)$, their external product is
\[
[X_1 , f_1]\times [X_2 , f_2 ] = [X_1 \times X_2 ,f_1 \times f_2]
\in \mathcal U^{d_1 + d_2} (Y_1\times Y_2).
\]
If admissible complex orientable Fredholm map $f_1 \times f_2$ is
transverse to the diagonal imbedding $\Delta:Y\rightarrow Y\times
Y$, then we do have an internal (cup) product
\[
[X_1 , f_1] \cup [X_2 , f_2] = \Delta^* [X_1\times X_2, f_1\times
f_2].
\]
If the diagonal imbedding $\Delta:Y\rightarrow Y\times Y$ is not
transverse to smooth proper Fredholm map $f_1 \times f_2:X_1
\times X_2 \rightarrow $, by Quinn's transversality Theorem, we
can find a smooth $\varepsilon$-approximation $\Delta'$ of
$\Delta$ which is transverse to $f_1\times f_2$. Then
\begin{thm}
If $[X_1,f_1] \in \mathcal U^{d_1}(Y_1)$ and $[X_2 , f_2]\in
\mathcal U^{d_2}(Y_2)$,internal(cup) product
\[
[X_1 , f_1] \cup [X_2 , f_2] = \Delta^* [X_1\times X_2, f_1\times
f_2] = \Delta'^* [X_1\times X_2, f_1\times f_2]\in \mathcal U^{d_1
+ d_2} (Y)
\]
where $\Delta'$ is a smooth $\varepsilon$-approximation of
$\Delta$ which is transverse to $f_1\times f_2$.
\end{thm}
The cup product is well-defined and associative.

Then, $\mathcal U^* (\,)$ is a multiplicative contravariant
functor for smooth functions on the separable Hilbert manifolds.

 We define the Euler class of a finite dimensional complex
vector bundle on a separable Hilbert manifold. Note that Theorem
\ref{son} implies that this Euler class is a well-defined
invariant of the bundle $\pi$.
\begin{defn}
Let $\pi: \xi \rightarrow B$ be a finite dimensional complex
vector bundle of dimension $d$ on a separable Hilbert manifold $B$
with zero-section $i: B\rightarrow \xi$. The \emph{$\mathcal
U$-theory Euler class} of $\xi$ is the element
$$
\chi(\pi) =i^* i_* (1)\in \mathcal U^{2d}(B).
$$
\end{defn}

Let $\pi: \xi \rightarrow X$ be a finite dimensional complex
vector bundle of dimension $d$ on a separable Hilbert manifold $X$
with zero-section $i: X\rightarrow \xi$.

Now we need a useful lemma from \cite{Quinn}.
\begin{lem}\label{tub}
A smooth split submanifold of a smooth separable Hilbert manifold
has a smooth tubular neighborhood.
\end{lem}
The map $i$ is proper so
that we have the Gysin map
$$
i_*: \mathcal U^{j}(X)\rightarrow \mathcal U^{j+2d}(\xi,\xi-U)
$$
where $U$ is a smooth neighborhood of the zero section.

The map $\pi$ is not proper. However if $U$ is contained in a tube
$U^r$ of finite radius $r$, then $\pi_{|\bar{U}}$ is proper and we
can define
$$
\pi_*:\mathcal U^{j+2d}(\xi,\xi-U)\rightarrow \mathcal U^{j}(X).
$$
Since $\pi i= \Id$ we have $\pi_* i_* = \Id$. The composite map
$i\pi $ is homotopic to $\Id_{\xi}$. If $U = U^{\circ}$ is itself
a tube, the homotopy moves on $U$ and we have \emph{Thom
isomorphism}
$$
\mathcal U^{j+2d}(\xi,\xi-U)\cong \mathcal U^{j}(X).
$$

\section{The ring structure on $LM^{-TM}$ and the Chas-Sullivan loop product in $U^*(LM)$.}

Let $M_d$ be a closed oriented $d$-dimensional smooth manifold,
and let $LM = C^{\infty}(\S^1 , M)$ be the space of smooth loops
in $M$. In \cite{chas}, Chas and Sullivan described an
intersection product on the homology $H_* (LM)$, having total
degree $-d$,
$$
\circ : H_q (LM) \otimes H_r (LM) \rightarrow H_{q + r -d} (LM).
$$

In \cite{cohen}, R. Cohen and J. Jones described a realization of
the Chas-Sullivan loop product in terms of a ring spectrum
structure on the Thom spectrum of a certain virtual bundle over
the loop space. We want to extend this product on the $U^*$-
theory.

Let $M_d$ be a closed complex $d$-dimensional smooth manifold, and
let $LM = C^{\infty}(\S^1 , M)$ be the space of smooth loops in
$M$. Let consider the standard parameterization of the circle by
the unit interval, $\exp: [0,1]\rightarrow \S^1$ defined by $\exp
(t) = \e^{2\pi it}$. With respect to this parameterization we can
regard a loop $\gamma \in LM$ as a map $\gamma: [0,1]\rightarrow
M$ with $\gamma(0) = \gamma(1)$. Let consider the evaluation map
$\text{ev}:LM \rightarrow M$ by $\gamma \rightarrow \gamma(1)$.

Let $\iota : M\rightarrow \C^{N+d}$ be a fixed smooth imbedding of
$M$ into codimension $N$ Unitary space. Let $\nu^{N} \rightarrow
M$ be the $2N$-dimensional normal bundle. Let $\text{Th}(\nu^N)$
be the Thom space of this bundle. We know that $\text{Th}(\nu^N)$
is Spanier-Whitehead dual to $M_+$ where $M_+$ denotes $M$ with a
disjoint basepoint. Let $M^{-TM}$ be the spectrum given by
desuspending the Thom space,
$$
M^{-TM}= \sum^{-2(N+d)}\text{Th}(\nu^N).
$$
We have the following spectra maps
$$
\S^0\rightarrow M_+ \wedge M^{-TM} \quad \text{and} \quad M_+
\wedge M^{-TM}\rightarrow \S^0,
$$
where $M^{-TM}$ is $\S$-dual of $M_+$. These maps induce an
equivalence with the function spectrum $M^{-TM} \simeq \text{Map}
(M_+,\S^0)$. Since $U^* (X) \simeq MU^*(X)$ for finite dimensional
manifolds $X$ and the contravariant cobordism theory $MU^*$ is
dual to the covariant bordism theory $MU_*$, we have the following
isomorphisms
\begin{align}
U^q(M_+) &\cong U_{-q}(M^{-TM})\notag\\
U^{-q}(M^{-TM}) &\cong U_q(M_+) \notag
\end{align}
for all $q\in \Z$. These duality isomorphisms are induced by the
compositions
$$
U^{-q}(M^{-TM}) \stackrel{\tau}\rightarrow
U^{-q+2d}(M_+)\stackrel{\rho}\rightarrow U_q(M_+)
$$
where $\tau$ is the Thom isomorphism, and $\rho$ is the Poincar`e
duality isomorphism for compact manifolds.

By duality, the diagonal map $\Delta: M\rightarrow M\times M$
induces a map of spectra
$$
\Delta^*: M^{-TM}\wedge M^{-TM} \rightarrow M^{-TM}
$$
that makes $ M^{-TM}$ into a ring spectrum with unit
$\S^0\rightarrow M^{-TM}$.
\par

Let $\text{Th}(\text{ev}^*(\nu^N))$ be the Thom space of the pull
back bundle $\text{ev}^*(\nu^N)\rightarrow LM$ where $LM =
C^{\infty}(\S^1 , M)$ is smooth manifold over
separable Hilbert space $\H$. Let define the spectrum
$$
LM^{-TM}= \sum^{-2(N+d)}\text{Th}(\text{ev}^*(\nu^N)).
$$

The representing dual manifold $X$ of the spectrum $LM^{-TM}$ is
also a smooth separable Hilbert manifold, e.g.$LM^{-TM}=[X,f]\in
U^{2d}(LM)$.

Now we will give the main theorem of this work.

\begin{thm}
The spectrum $LM^{-TM}$ is a homotopy commutative ring spectrum
with unit,whose multiplication
$$
\mu:LM^{-TM}\wedge LM^{-TM}\rightarrow LM^{-TM}
$$
satisfies the following properties.
\par
\bf{1.} The evaluation map $\text{ev}:LM^{-TM}\rightarrow M^{-TM}$
is a map of ring spectra.

\par
\bf{2.} There is a map of ring spectra $\rho:LM^{-TM}\rightarrow
\ds\sum^{\infty}\Omega M$ where the target is the suspension
spectrum of the based loop space with a disjoint basepoint. Its
ring structure is induced by the usual product on the based loop
space. In bordism the map $\rho_*$ is given by the composition
$$
\rho_*: MU_q (LM^{-TM})\stackrel{\tau}\rightarrow
MU_{q+2d}(LM)\stackrel{\iota}\rightarrow MU_q(\Omega M)
$$
where $\tau$ is the Thom isomorphism and $\iota$ takes a bordism
class with dimension $(q+2d)$ and by intersecting with the based
loop $\Omega M$ as a codimension $2d$,e.g. $\iota=
i^*:MU_*(LM)\rightarrow MU_{*-2d}(\Omega M)$ is an induced
homomorphism from the embedding $i:\Omega M\rightarrow LM$.

\par
\bf{3.} The ring structure is compatible with the Chas-Sullivan
loop product in the sense that the following diagrams commute.
$$
\begin{CD}
U^q (LM^{-TM})\times U^r (LM^{-TM})@>\textbf{ext}>> U^{q +r}
(LM^{-TM}\wedge LM^{-TM})  @>\Delta^*>>  U^{q+r} (LM^{-TM}) \\
@V\cong V{u_*}V   @VVV      @V\cong V{u_*}V \\   U^{q
+2d}(LM)\times U^{r+2d} (LM) @>>{\circ}>\ldots @>>{\circ}>
U^{q+r+2d}(LM)
\end{CD}
$$
and
$$
\begin{CD}
U_{q} (LM^{-TM})\times U_r (LM^{-TM})@>\textbf{ext}>> U_{q +r}
(LM^{-TM}\wedge LM^{-TM})  @>\mu_*>>  U_{q+r} (LM^{-TM}) \\
@V\cong V{u_*}V   @VVV      @V\cong V{u_*}V \\   U_{q
+2d}(LM)\times U_{r+2d} (LM) @>>{\circ}>\ldots @>>{\circ}>
U_{q+r+2d}(LM)
\end{CD}
$$

where $\textbf{ext}$ is the external product, $u_{*}$ is the Thom
isomorphism, $\Delta:LM^{-TM}\times LM^{-TM}\rightarrow LM^{-TM}$
is the diagonal map which is adjoint to the multiplication map
$$\mu:LM^{-TM}\wedge LM^{-TM}\rightarrow LM^{-TM}$$ and $\circ$ is
the Chas-Sullivan loop product in cobordism.

\end{thm}
\begin{proof}
The proof was done in \cite{cohen} by essentially the same
argument for homology but their proof had a fundamental mistake.
In this proof we will sort out this mistake and we will do the
modification for cobordism and bordism theories.
\par
Let $\Delta:M\rightarrow M\times M$ be the diagonal embedding of
closed oriented manifold $M$. The normal bundle is isomorphic to
the tangent bundle, $\nu_{\Delta}\cong TM$ so that the Pontrjagin-
Thom map is a complex orientable map $\tau: M\times M\rightarrow
M^{TM}$ with index zero. So we have Gysin map in cobordism,
$$
MU^* (M\times M)\stackrel{\tau_*}\rightarrow
MU^*(M^{TM})\stackrel{u_*}\rightarrow MU^{*-2d}(M)
$$
which is the transversal intersection product.

\par
Here we will apply the Pontrjagin-Thom construction to the
diagonal embedding $\Delta:M\rightarrow M\times M$ using the
canonical bundle $-TM \times - TM$ over $M\times M$. We get a map
of Thom spectra
$$
\tau:(M\times M)^{-TM \times - TM} \rightarrow M^{TM \oplus
\Delta^*(-TM \times - TM)}
$$
or,
$$
\tau: M^{-TM} \wedge M^{-TM} \rightarrow M^{-TM}.
$$
The details about the Pontrjagin-Thom construction can be found in
\cite{cohen}.

To construct the ring spectrum product
$$
\mu:LM^{-TM}\wedge LM^{-TM}\rightarrow LM^{-TM},
$$
they pull back the structure $\tau$ over the loop space $LM$.

For this, they define $LM\times_M LM$ which is fiber product in
the following diagram
$$
\begin{CD}
LM\times_M LM  @>\bar{\Delta}>> LM\times LM  \\
@V\textbf{ev}VV    @VV\textbf{ev}\times\textbf{ev} V \\
M @>>\Delta > M\times M.
\end{CD}
$$
They note that $LM\times_M LM$ is a codimension $2d$ submanifold
of the infinite dimensional manifold $LM\times LM$ and it is equal
to
$$
\{(\alpha,\beta)\in LM\times LM  : \alpha(0)=\beta(0)\}.
$$
Since $\textbf{ev}:LM\rightarrow M$ is a submersion, the fiber
product corresponds to transversal intersection of maps, so
$\bar{\Delta}$ is pull back of the diagonal map $\Delta$ under the
submersion map $\textbf{ev}\times\textbf{ev}: LM\times
LM\rightarrow M\times M$. The induced map $\bar{\Delta}$ is a
Fredholm map with index $2d$ and consequently $LM\times_M LM$ is a
codimension $2d$ smooth submanifold of the infinite dimensional
manifold $LM\times LM$.

Also they note that there is a natural map $\gamma:LM\times_M LM
\rightarrow LM$ defined by first applying $\alpha$ and then
$\beta$. That is,
$$
\gamma:LM\times_M LM \rightarrow LM \quad \text{by} \quad
\gamma((\alpha,\beta))=\alpha\ast\beta,
$$
where
$$
(\alpha\ast\beta) (t) =\begin{cases} \alpha(2t) & \text{if}\quad
0\leqq t\leqq \frac{1}{2}\\
\beta(2t-1) & \text{if}\quad \frac{1}{2}\leqq t\leqq 1.
\end{cases}
$$
But the image of $\gamma$ is not smooth loop as they have
described it above. However there is a standard way to modify the
definition of $\gamma$ so that the target is smooth. The
resolution is in the parametrization of the loop. With continuous
loops one just defines:
$$
(\alpha\ast\beta) (t) =\begin{cases} \alpha(2t) & \text{if}\quad
0\leqq t\leqq \frac{1}{2}\\
\beta(2t-1) & \text{if}\quad \frac{1}{2}\leqq t\leqq 1.
\end{cases}
$$
but this may be not be smooth at $\frac{1}{2}$ or at $1$. To
modify this we need a bijective smooth function from $[0,1]$ to
$[0,1]$ which has all derivatives zero at $0$ and $1$ since we
will use this to reparameterize the two loops. Then we can patch
them together without losing smoothness.

To write this special parametrization, we define a map $\varphi:
[0,1]\rightarrow [0,1]$ by
$$
\varphi(t)=\frac{1}{c}{\ds\int_{0}^{t}\exp^{(-\frac{1}{s^2(s-1)^2})}ds}
\quad \text{where}\quad
c={\ds\int_{0}^{1}\exp^{(-\frac{1}{s^2(s-1)^2})}ds}.
$$
It is a smooth bijective function from $[0,1]$ to $[0,1]$ which
has all derivatives zero at $0$ and $1$. Now we can define a map
$$
\Phi:LM\rightarrow LM \quad \text{by}\quad
\Phi(\alpha)(t)=\alpha(\varphi(t))=\alpha\circ\varphi(t).
$$
It is a smooth map and let $L_{\rs}M =\Phi(LM)$. $L_{\rs}M$ is a
smooth submanifold of $LM$ and it is smooth homotopic retraction
of $LM$. Similarly $L_{\rs}M\times_{M}L_{\rs}M$ can be
constructed. It is also smooth homotopic retraction of $LM\times_M
LM$. Then we can define a new version of the map
$\gamma:LM\times_M LM\rightarrow LM$ as the composition of the
following maps
$$
\gamma:LM\times_M LM\stackrel{\Phi\times\Phi}\rightarrow
L_{\rs}M\times_{M}L_{\rs}M\stackrel{\ast}\rightarrow
L_{\rs}M\subset LM.
$$
The defined new version of $\gamma$ is smooth. In the proof of
\cite{cohen}, we can use $L_{\rs}M\times_{M}L_{\rs}M$ instead of
$LM\times_M LM$ because it is smooth homotopic retraction of
$LM\times LM$.

If we restrict to the product of the based loop spaces,
$\Omega_{\rs}M \times \Omega_{\rs}M \subset
L_{\rs}M\times_{M}L_{\rs}M$, then $\gamma$ is just the $H$-space
product on the based loop space, $\Omega_{\rs}M \times
\Omega_{\rs}M\rightarrow \Omega_{\rs}M$.

\par

The embedding
$\bar{\Delta}:L_{\rs}M\times_{M}L_{\rs}M\hookrightarrow
L_{\rs}M\times L_{\rs}M$ has a tubular neighborhood
$\nu(\bar{\Delta})$ defined to be the inverse image of the tubular
neighborhood of the diagonal map $\Delta: M\hookrightarrow M\times
M$:
$$
\nu(\bar{\Delta})=\textbf{ev}^{-1}(\nu(\Delta)).
$$
Hence there is a Pontjagin-Thom construction
$$
\tau:L_{\rs}M\times L_{\rs}M \rightarrow
(L_{\rs}M\times_{M}L_{\rs}M)^{\textbf{ev}^*(TM)}.
$$
The map $\tau$ is a smooth Fredholm map. By the Pontrjagin-Thom
construction, we have the following commutative diagram
$$
\begin{CD}
L_{\rs}M\times L_{\rs}M  @>{\tau}>> (L_{\rs}M\times_M L_{\rs}M)^{TM}  \\
@V\textbf{ev}VV    @VV\textbf{ev} V \\
M\times M @>>\tau > M^{TM}.
\end{CD}
$$
In bordism, we have
$$
\iota:U_*(L_{\rs}\times L_{\rs}M)\rightarrow
U_{*-2d}(L_{\rs}\times_M L_{\rs}M)
$$
where $\iota$ takes a bordism class with dimension $n$ and
intersects with the submanifold $L_{\rs}\times_M L_{\rs}M$ as a
codimension $2d$, i.e. pull backs by the inclusion
$L_{\rs}\times_M L_{\rs}M \rightarrow L_{\rs}\times L_{\rs}M$.

\par
By the following commutative diagram
$$
\begin{CD}
L_{\rs}M\times_M L_{\rs}M  @>{\gamma}>> L_{\rs}M  \\
@V\textbf{ev}VV    @VV\textbf{ev} V \\
M @>>{=} > M,
\end{CD}
$$
we have an induced map of bundles
$\gamma:\textbf{ev}^*(TM)\rightarrow \textbf{ev}^*(TM)$,
hence we have a map of spectra
$$
(L_{\rs}M\times_M L_{\rs}M)^{TM}  \stackrel{\gamma}\rightarrow (L_{\rs}M)^{TM}.
$$

Then we will get the following composition
$$
\tilde{\mu}:L_{\rs}M\times L_{\rs}M  \stackrel{\tau}\rightarrow
(L_{\rs}M\times_M L_{\rs}M)^{TM}\stackrel{\gamma}\rightarrow LM^{TM}.
$$
In bordism, the homomorphism
$$
U_{*}(LM\times LM)\stackrel{\tilde{\mu}}\rightarrow U_*(LM^{TM})
\stackrel{u_*}\rightarrow U_{*-2d}(LM)
$$
takes a bordism class in $LM \times LM$, intersects in with the codimension d submanifold
$L_{\rs}M\times_M L_{\rs}M$, maps it via $\gamma$ to $LM$. This is the definition of
Chas-Sullivan product $U_*(LM)$.

\par

Using the diagonal embedding $LM\rightarrow LM\times LM$, we can
perform the Pontrjagin-Thom construction when we pull back the
virtual bundle $-TM\times -TM$ over $LM\times LM$. Then we obtain
$$
\tau: LM^{-TM}\wedge LM^{-TM} \rightarrow (L_{\rs}M\times_M
L_{\rs}M)^{TM \oplus -2TM} =(L_{\rs}M\times_M L_{\rs}M)^{-TM}.
$$
Then we can define the ring structure on the Thom  spectrum to be
the composition
$$
\mu: LM^{-TM}\wedge LM^{-TM} \stackrel{\tau}\rightarrow
(L_{\rs}M\times_M L_{\rs}M)^{-TM}\stackrel{\gamma}\rightarrow
LM^{-TM}.
$$
In \cite{cohen}, They show that $\mu$ is associative.
\par
In bordism, by Thom isomorphism $\mu_{*}$ induces the same
homomorphism as $\tilde{\mu}_*$, so we have the following diagram
commutes.
$$
\begin{CD}
U_{q-4d}(LM^{-TM}\wedge LM^{-TM})  @>{\mu_*}>> U_{q-4d}(LM^{-TM})  \\
@V{u_*}V{\cong}V    @V{u_*}V{\cong} V \\
U_{q}(LM \times LM) @>>{\circ} >U_{q-2d}(LM),
\end{CD}
$$
where $\circ:U_{q}(LM \times LM) \rightarrow U_{q-2d}(LM)$ is the
Chas-Sullivan product. In complex cobordism, we define the
Chas-Sullivan product by the following commutative diagram
$$
\begin{CD}
U^{q-4d}(LM^{-TM}\wedge LM^{-TM})  @>{\Delta^*}>> U^{q-4d}(LM^{-TM})  \\
@V{u_*}V{\cong}V    @V{u_*}V{\cong} V \\
U^{q}(LM \times LM) @>>{\circ} >U^{q-2d}(LM),
\end{CD}
$$
In \cite{cohen}, They show that $\rho:LM^{-TM}\rightarrow
\sum^{\infty}(\Omega M))$ is a map of ring spectra. In bordism the
map $\rho_*$ is given by the composition
$$
\rho_*: MU_q (LM^{-TM})\stackrel{\tau}\rightarrow
MU_{q+2d}(LM)\stackrel{\iota}\rightarrow MU_q(\Omega M)
$$
where $\tau$ is the Thom isomorphism and $\iota$ takes a bordism
class with dimension $(q+2d)$ and by intersecting with the based
loop $\Omega M$ as a codimension $2d$,e.g. $\iota=
i^*:MU_*(LM)\rightarrow MU_{*-2d}(\Omega M)$ is an induced
homomorphism from the embedding $i:\Omega M\rightarrow LM$.

\end{proof}

\end{document}